# A CLASS OF STATIONARY SEQUENCES

Florentin Smarandache, Ph D
Department of Math & Sciences
University of New Mexico
200 College Road
Gallup, NM 87301, USA
E-mail: smarand@unm.edu

§1. We define a class of sequences $\{a_n\}$ by $a_1 = a$ and $a_{n+1} = P(a_n)$, where $P$ is a polynomial with real coefficients. For which $a$ values, and for which polynomials P will these sequences be constant after a certain rank? Then we generalize it from polynomials P to real functions f.
  In this note, the author answers this question using as reference F. Lazebnik & Y. Pilipenko's E 3036 problem from A. M. M., Vol. 91, No. 2/1984, p. 140.
  An interesting property of functions admitting fixed points is obtained.

§2. Because $\{a_n\}$ is constant after a certain rank, it results that $\{a_n\}$ converges. Hence, $(\exists) e \in R: e = P(e)$, that is the equation $P(x) - x = 0$ admits real solutions. Or $P$ admits fixed points $((\exists) x \in R: P(x) = x)$.
  Let $e_1, ..., e_m$ be all real solutions of this equation. We construct the recurrent set $E$ as follows:
  1) $e_1, ..., e_m \in E$;
  2) if $b \in E$ then all real solutions of the equation $P(x) = b$ belong to $E$;
  3) no other element belongs to $E$, except those elements obtained from the rules 1) and/or 2), applied for a finite number of times.
  We prove that this set $E$, and the set $A$ of the "$a$" values for which $\{a_n\}$ becomes constant after a certain rank, are indistinct.

Let's show that "$E \subseteq A$":
  1) If $a = e_i$, $1 \le i \le m$, then $(\forall) n \in \mathbb{N}^*$ $a_n = e_i = $ constant.
  2) If for $a = b$ the sequence $a_1 = b$, $a_2 = P(b)$ becomes constant after a certain rank, let $x_0$ be a real solution of the equation $P(x) - b = 0$, the new formed sequence: $a'_1 = x_0$, $a'_2 = P(x_0) = b$, $a'_3 = P(b), ...$ is indistinct after a certain rank with the first one, hence it becomes constant too, having the same limit.
  3) Beginning from a certain rank, all these sequences converge towards the same limit $e$ (that is: they have the same $e$ value from a certain rank) are indistinct, equal to $e$.

Let's show that "$A \subseteq E$":

Let "$a$" be a value such that: $\{a_n\}$ becomes constant (after a certain rank) equal to $e$. Of course $e \in \{e_1, ..., e_m\}$ because $e_1, ..., e_m$ are the only values towards these sequences can tend.

If $a \in \{e_1, ..., e_m\}$, then $a \in E$.

Let $a \notin \{e_1, ..., e_m\}$, then $(\exists) n_0 \in \mathbb{N}: a_{n_0+1} = P(a_{n_0}) = e$, hence we obtain by applying the rules 1) or 2) a finite number of times. Therefore, because $e \in \{e_1, ..., e_m\}$ and the equation $P(x) = e$ admits real solutions we find $a_{n_0}$ among the real solutions of this equation: knowing $a_{n_0}$ we find $a_{n_0-1}$ because the equation $P(a_{n_0-1}) = a_{n_0}$ admits real solutions (because $a_{n_0} \in E$ and our method goes on until we find $a_1 = a$ hence $a \in E$.

**Remark.** For $P(x) = x^2 - 2$ we obtain the E 3036 Problem (A. M. M.).
Here, the set $E$ becomes equal to

$$\{\pm 1, 0, \pm 2\} \cup \left\{ \pm\underbrace{\sqrt{2 \pm \sqrt{2 \pm ... \sqrt{2}}}}_{n_0 \text{ times}}, n \in \mathbb{N}^* \right\} \cup \left\{ \pm\underbrace{\sqrt{2 \pm \sqrt{2 \pm ... \sqrt{2 \pm \sqrt{3}}}}}_{n_0 \text{ times}}, n \in \mathbb{N} \right\}$$

Hence, for all $a \in E$ the sequence $a_1 = a$, $a_{n+1} = a_n^2 - 2$ becomes constant after a certain rank, and it converges (of course) towards $-1$ or $2$:

$$(\exists) n_0 \in \mathbb{N}^*: (\forall) n \geq n_0 \quad a_n = -1$$

or

$$(\exists) n_0 \in \mathbb{N}^*: (\forall) n \geq n_0 \quad a_n = 2.$$

**Generalization**.

This can be generalized to defining a class of sequences $\{a_n\}$ by $a_1 = a$ and $a_{n+1} = f(a_n)$, where f: R $\to$ R is a real function. For which $a$ values, and for which functions f will these sequences be constant after a certain rank?

In a similar way, because $\{a_n\}$ is constant after a certain rank, it results that $\{a_n\}$ converges. Hence, $(\exists) e \in R: e = f(e)$, that is the equation $f(x) - x = 0$ admits real solutions. Or f admits fixed points $((\exists) x \in R: f(x) = x)$.

Let $e_1, ..., e_m$ be all real solutions of this equation. We construct the recurrent set $E$ as follows:
  1) $e_1, ..., e_m \in E$;
  2) if $b \in E$ then all real solutions of the equation $f(x) = b$ belong to $E$;
  3) no other element belongs to $E$, except those elements obtained from the rules 1) and/or 2), applied for a finite number of times.

Analogously, this set $E$, and the set $A$ of the "$a$" values for which $\{a_n\}$ becomes constant after a certain rank, are indistinct.